\documentclass[12pt]{amsart}
\usepackage{amsfonts}
\usepackage{amsmath}
\usepackage{amssymb}
\usepackage[mathcal]{eucal}

\usepackage[bookmarks,colorlinks,breaklinks]{hyperref}

\newtheorem{theorem}{Theorem}[section]
\newtheorem{corollary}[theorem]{Corollary}
\newtheorem{lemma}[theorem]{Lemma}
\newtheorem{proposition}[theorem]{Proposition}
\newtheorem{thmA}{Theorem}

\theoremstyle{definition}
\newtheorem{question}[theorem]{Question}
\newtheorem{definition}[theorem]{Definition}
\newtheorem{case}{Case}
\numberwithin{equation}{section}
\theoremstyle{remark}
\newtheorem{remark}[theorem]{Remark}
\def\Aut{\operatorname{Aut}}
\def\Sym{\operatorname{Sym}}
\def\hull#1{\langle#1\rangle}
\def\cd{\mathrm{cd}}
\def\dd{\mathrm{dd}}
\def\JM{\mathbf{JM}}
\def\ad{\operatorname{ad}}
\def\cof{\operatorname{cf}}
\begin{document}

	\title[Density Spectra]{Density Spectra of Topological Groups}
	
\author[D. Peng]{Dekui Peng\textsuperscript{*}}
 \address[D. Peng]
{Institute of Mathematics, Nanjing Normal University, Nanjing 210024, China}
\email{pengdk10@lzu.edu.cn}
\thanks{*Corresponding author.}

\author[Y. Zhou]{Yi Zhou}
 \address[Y. Zhou]
{School of Sciences, Southwest Petroleum University,
Chengdu 610500, China}
\email{zhouyi.math@outlook.com}

	
\begin{abstract}
We study density spectra associated with dense and closed subgroups of
topological groups. For a topological group $G$, let $\dd^*(G)$ denote the
set of densities of its dense subgroups. We prove that every compact group
satisfies
$\dd^*(G)=[d(G),w(G)]$. 

We also investigate the density spectrum $\cd(G)$ of infinite closed
subgroups. In $\mathbf{ZFC}$, we construct a separable countably compact
Boolean group containing a closed non-separable subgroup, answering a
question of Leiderman, Morris, and Tkachenko. We further show that a free
pro-$p$ group of uncountable rank has no non-trivial metrizable closed
normal subgroup, giving a negative answer to a problem of Hern\'andez,
Hofmann, and Morris. Finally, we study when $\cd(G)$ is an interval. Among other
results, we show that
under Shelah's Strong Hypothesis ($\mathbf{SSH}$), the closed density spectrum of every
infinite $\omega$-bounded group of countable tightness is an interval of cardinals.
\end{abstract}
\maketitle
\setcounter{tocdepth}{1}
\section{Introduction}\label{s1}

The density character, or simply the density, of a topological space $X$,
denoted by $d(X)$, is one of the basic cardinal invariants measuring its
topological size. In contrast with weight, density behaves rather subtly
with respect to subspaces. In particular, neither dense nor closed
subspaces need have density comparable in an obvious way with that of the
ambient space. The familiar failure of hereditary separability already
illustrates this phenomenon: a separable space may contain a closed
non-separable subspace.

This naturally leads to the problem of determining which densities can
occur among distinguished subspaces of a fixed space. Juh\'asz, van Mill,
Soukup, and Szentmikl\'ossy introduced in \cite{JvMSS} the
\emph{double density spectrum}
\[
\dd(X)=\{d(Y):Y\subseteq X\text{ is dense}\}.
\]
More recently, Juh\'asz and van Mill proved in \cite{JvM} that
$\dd(X)=[d(X),w(X)]$ whenever $X$ is homeomorphic to a locally compact
group. Thus, from the point of view of arbitrary dense subspaces, locally
compact groups exhibit the largest possible density spectrum.

For a topological group, however, it is natural to require the dense
subspaces under consideration to respect the algebraic structure. We
therefore consider the following refinement.

\begin{definition}
Let $G$ be a topological group. The \emph{density spectrum of dense
subgroups} of $G$ is
\[
\dd^*(G)=\{d(H):H\leq G\text{ is dense}\}.
\]
\end{definition}

Following the terminology motivated by \cite{JvM}, we use the following
notation.

\begin{definition}
A Tychonoff space $X$ is said to have property $\JM$ if
$\dd(X)=[d(X),w(X)]$. A topological group $G$ is said to have property
$\JM^*$ if
$\dd^*(G)=[d(G),w(G)]$.
\end{definition}

The passage from $\JM$ to $\JM^*$ is substantial. Indeed, although every
locally compact group satisfies $\JM$, the corresponding statement for
dense subgroups is false even within natural subclasses of locally compact
groups. It was proved in \cite{Peng} that there are non-metrizable locally
compact groups, including abelian and connected examples, whose every
dense subgroup is separable. Such groups fail property $\JM^*$. On the
positive side, the same paper showed that $w(G)\in\dd^*(G)$ when $G$ is
compact abelian, compact connected, or connected of weight greater than
$2^\omega$. Related questions concerning possible densities and weights
of quotient groups of precompact abelian groups were studied in
\cite{Peng1}.

These results leave open a natural compact-group problem: does every
compact group satisfy $\JM^*$? Our first main theorem gives a complete
affirmative answer.

\begin{thmA}\label{ThA}
Every compact group satisfies property $\JM^*$. Equivalently, for every
compact group $G$,
\[
\dd^*(G)=[d(G),w(G)].
\]
\end{thmA}

The second part of the paper concerns a complementary problem. For an
infinite topological group $G$, define the \emph{density spectrum of closed
subgroups} by
\[
\cd(G)=\{d(H):H\leq G\text{ is an infinite closed subgroup}\}.
\]
Dense and closed subgroups display markedly different behavior. A
classical theorem of Comfort and Itzkowitz \cite{CIt} implies, in
particular, that closed subgroups of separable locally compact groups are
separable. Beyond local compactness this phenomenon breaks down sharply.
Leiderman, Morris, and Tkachenko \cite{LMT} showed that closed subgroups of
separable precompact, and even pseudocompact, groups may have arbitrarily
large density within the natural cardinal restrictions. They also asked
whether a separable countably compact group can contain a closed
non-separable subgroup in $\mathbf{ZFC}$.

Using the recent construction of Hru\v{s}\'ak, van Mill, Ramos-Garc\'ia,
and Shelah \cite{HMRS} of a countably compact Boolean group without
non-trivial convergent sequences, we answer this question affirmatively.

\begin{thmA}\label{ThB}
In $\mathbf{ZFC}$, there exists a separable countably compact Boolean group
containing a closed non-separable subgroup.
\end{thmA}

Our construction combines the group supplied by \cite{HMRS} with a
suitable $\Sigma$-product inside a Boolean cube. The resulting group
remains countably compact and separable, while the $\Sigma$-product
appears as a closed subgroup of density $\mathfrak c$. Thus the contrast
between the ambient group and its closed subgroup occurs already for
Boolean groups and requires no additional set-theoretic assumption.

We next return to compact groups. Hern\'andez, Hofmann, and Morris
\cite{HHM} proved that every infinite locally compact group $G$ has closed
subgroups of every prescribed weight between $\omega$ and $w(G)$. For a
compact group, using $d(H)=\log w(H)$ for infinite closed subgroups $H$,
their theorem gives the exact description
\[
\cd(G)=\{\log\kappa:\kappa\in[\omega,w(G)]\}.
\]
They also asked whether the closed subgroup of prescribed weight can
always be chosen normal. We show that this is false in a strong form:
free pro-$p$ groups of uncountable rank have no non-trivial metrizable
closed normal subgroups. The proof is an application of the normal-subgroup theorem for free pro-$p$ groups.

The description above also shows that, even for compact groups, $\cd(G)$
need not be an interval of cardinals. This raises the question of how much
interval structure can nevertheless be forced by additional topological
assumptions. We obtain first a general result in terms of tightness. If
$G$ is a topological group of density $\kappa$ and $\tau$ is regular with
$t(G)<\tau\leq\kappa$, then $G$ contains a closed subgroup of density
$\tau$. Consequently, for a group of countable tightness, every regular
cardinal below a density already occurring in $\cd(G)$ also belongs to
$\cd(G)$.

The singular case is more delicate. For an $\omega$-bounded group $H$ of
countable tightness and density $\lambda$, we prove the estimate
\[
d(L)\leq
\max\left\{\omega,
\operatorname{cf}\bigl([\lambda]^{\leq\omega},\subseteq\bigr)\right\}
\]
for every closed subgroup $L$ of $H$. The point is that the closure of
every countably generated subgroup is compact and, by countable
tightness, metrizable. This reduces the singular-cardinal problem to the
cofinal structure of the poset $[\lambda]^{\leq\omega}$.

As a consequence, if $\tau$ is singular and
$\operatorname{cf}([\lambda]^{\leq\omega},\subseteq)<\tau$ for every
infinite $\lambda<\tau$, then an $\omega$-bounded group of countable
tightness cannot have a gap at $\tau$: whenever
$\kappa\in\cd(G)$ and $\kappa\geq\tau$, one also has
$\tau\in\cd(G)$. This applies, in particular, when
$\lambda^\omega<\tau$ for every $\lambda<\tau$, and hence to all
singular strong limit cardinals. It also yields in $\mathbf{ZFC}$ the
notable conclusion that $\aleph_\omega$ cannot be such a gap.

The same criterion connects the problem naturally with pcf theory.
Recall that Shelah's Strong Hypothesis $\mathbf{SSH}$ asserts
$pp(\mu)=\mu^+$ for every singular cardinal $\mu$. 

\begin{thmA}\label{ThC}
Assume $\mathbf{SSH}$. Then for every infinite $\omega$-bounded topological group
$G$ with countable tightness, $\cd(G)$ is an interval of cardinals.
\end{thmA}

Thus the remaining obstruction in $\mathbf{ZFC}$ for $\omega$-bounded
groups is concentrated at singular cardinals with large countable-covering
cofinalities. For general infinite countably compact groups, even the
compact-metrizable hull argument used above is unavailable. We conclude
by asking whether $\cd(G)$ must be an interval for infinite countably compact
groups of countable tightness, and whether the same conclusion holds in
$\mathbf{ZFC}$ for all infinite $\omega$-bounded groups of countable tightness.

Section~\ref{s3} develops these questions concerning closed subgroups,
including the $\mathbf{ZFC}$ countably compact construction, the
pro-$p$ counterexample, and the regular, singular, and pcf-theoretic
results described above.

\subsection{Notation and Terminology}

Throughout this paper, all topological spaces are assumed to be Tychonoff. 
The cardinality of a set $A$ is denoted by $|A|$. For cardinals $\lambda \le \tau$, the interval $[\lambda, \tau]$ denotes the set of all cardinals $\kappa$ such that $\lambda \le \kappa \le \tau$. It should be emphasized that all such intervals in this paper represent intervals of cardinals, rather than ordinals. For an infinite cardinal $\kappa$, its logarithm is defined as $\log\kappa = \min\{\lambda : 2^\lambda \ge \kappa\}$.  We write $\omega$ for the first infinite cardinal, i.e., $\omega=\{0,1,2,\ldots\}$.

The \emph{density character} (or simply the \emph{density}) of a topological space $X$, denoted by $d(X)$, is defined as the minimal cardinality of a dense subset of $X$:
\[ d(X) := \min \{ |Y| : Y \subseteq X \text{ is a dense subset} \}. \]
Another closely related and widely studied cardinal invariant is the \emph{topological weight}, $w(X)$, which is the minimal cardinality of a base for the topology of $X$. In the realm of compact groups, the weight $w(G)$ coincides with its character, i.e., the minimal cardinality of a local base at the identity. Moreover, if $G$ is an infinite compact group, then $d(G)=\log w(G)$ \cite{AT}.
For any infinite Tychonoff space $X$, it is a well-known fact that the inequalities $d(X) \le w(X) \le 2^{d(X)}$ hold \cite{Eng}.

For a subset $X$ of a topological group $G$, we denote the subgroup generated by $X$ by $\hull{X}$, and its topological closure is denoted by $\overline{\hull{X}}$. A topological group is called a \emph{Boolean group} if every non-identity element has order $2$. Clearly, Boolean groups are abelian.

Let $G$ be a topological abelian group. Its Pontryagin--van Kampen dual, denoted by $\widehat{G}$, is the group of all continuous characters from $G$ to the circle group $\mathbb{T}$, endowed with the compact-open topology. The word ``character'' has two meanings in this paper: earlier it denoted the cardinal invariant $\chi(X)$ measuring the size of local bases, whereas here and throughout the remainder of the paper it means a continuous homomorphism into $\mathbb{T}$. The Pontryagin--van Kampen Duality Theorem asserts that the natural evaluation map identifies $G$ with its double dual $\widehat{\widehat{G}}$ as topological groups when $G$ is locally compact.

If $f: G \to H$ is a continuous homomorphism between locally compact abelian groups, we denote by $\widehat{f}: \widehat{H} \to \widehat{G}$ the induced dual homomorphism, defined by $\widehat{f}(\chi) = \chi \circ f$. 
For a closed subgroup $H$ of $G$, its \emph{annihilator} $H^\perp$ is defined as
\[ H^\perp := \{ \chi \in \widehat{G} : \chi(h) = 0 \text{ for all } h \in H \}. \]
It is a standard fact that $\widehat{G/H}$ is topologically isomorphic to $H^\perp$. We shall use the standard properties of Pontryagin duality freely throughout this paper without further reference. For a comprehensive treatment, we refer the reader to standard monographs such as \cite{AT} or \cite{HR}.

A \emph{profinite group} is a compact, totally disconnected topological group. Equivalently, profinite groups can be characterized as inverse limits of inverse systems of finite groups equipped with the discrete topology. 
If $p$ is a prime and all finite groups in the inverse system are $p$-groups, then the inverse limit is called a \emph{pro-$p$ group}.
For a detailed and fruitful exposition on profinite groups, we refer to \cite{RZ}.

\subsection{Some Lemmas}

To prepare for the proof of the main theorems, we collect several auxiliary lemmas. These results describe how weights and densities behave under common group-theoretic and topological constructions. They will be used repeatedly in various reduction steps throughout the paper.

\medskip

Throughout, all product factors and semidirect product components are naturally identified with closed subgroups of the ambient group.

\medskip

We begin with a lemma concerning quotients by totally disconnected central subgroups in connected compact groups.

\begin{lemma}\label{conq}
Let $G$ be a connected compact group, and let $D$ be a closed, totally disconnected normal subgroup. Then
\[
w(G/D) = w(G).
\]
\end{lemma}

\begin{proof}
Since $G$ is connected and $D$ is totally disconnected and normal, $D$ is central in $G$. The asserted equality is a standard consequence of the structure theory of connected compact groups; see \cite[Proposition~12.19(iii)]{HM}.
\end{proof}

Next, we consider how the density of subgroups behaves under open continuous group homomorphisms.

\begin{lemma}\label{preimage}
Let $\pi: G \to H$ be a surjective open continuous homomorphism between locally compact groups.  
If $D$ is a dense subgroup of $H$ with density $\kappa \geq d(G)$, then the preimage $\pi^{-1}(D)$ is a dense subgroup of $G$ and satisfies
\[
d(\pi^{-1}(D)) = \kappa.
\]
\end{lemma}

\begin{proof}
Let $N = \ker(\pi)$. By a theorem of Comfort and Itzkowitz \cite[2.4 Theorem]{CIt}, we have $d(N) \leq d(G) \leq \kappa$.  
Using the standard inequality for the density of extensions,
\[
d(\pi^{-1}(D)) \leq d(D) \cdot d(N) \leq \kappa \cdot \kappa = \kappa.
\]
On the other hand, if $E$ is dense in $\pi^{-1}(D)$, then $\pi(E)$ is dense in $D$. Hence
\[
d(D)\leq d(\pi^{-1}(D)),
\]
so $d(\pi^{-1}(D))\geq\kappa$. Therefore $d(\pi^{-1}(D))=\kappa$.
\end{proof}

Finally, we recall a useful finite-support construction for building dense subgroups of prescribed density in products of compact metrizable groups.

\begin{definition}
Let $\{K_i : i \in I\}$ be a family of topological groups, and let $e_i \in K_i$ denote the identity element of $K_i$. The \emph{finite-support subgroup} of $\prod_{i\in I}K_i$ is
\[
\bigoplus_{i\in I}K_i := \left\{\, x \in \prod_{i \in I} K_i : \{i \in I : x_i \neq e_i\} \text{ is finite} \,\right\}.
\]
It is precisely the subgroup generated by the coordinate copies of the groups $K_i$.
\end{definition}

\begin{lemma}\label{prod}
Let $\{K_i : i \in I\}$ be an infinite family of non-trivial topological groups, such that $d(K_i)\leq |I|$ for each $i\in I$.
Then the finite-support subgroup $\bigoplus_{i\in I}K_i$ is dense in $\prod_{i \in I} K_i$ and has density equal to $|I|$.
\end{lemma}

\begin{proof}
Put $\kappa=|I|$. For every $i\in I$, choose a dense set $D_i\subseteq K_i$ such that $e_i\in D_i$ and $|D_i|\leq\kappa$. Let
\[
D=\left\{x\in\bigoplus_{i\in I}K_i:x_i\in D_i\text{ for every }i\in I\right\}.
\]
Then $|D|\leq\kappa$, and $D$ meets every non-empty basic open subset of $\prod_{i\in I}K_i$. Hence $\bigoplus_{i\in I}K_i$ is dense in the product and
\[
d\left(\bigoplus_{i\in I}K_i\right)\leq\kappa.
\]

Conversely, let $E$ be dense in $\bigoplus_{i\in I}K_i$. For each $i\in I$, choose a non-empty open set $U_i\subseteq K_i$ such that $e_i\notin U_i$. Since the cylinder $\pi_i^{-1}(U_i)$ meets $E$, there exists $x\in E$ with $i\in\operatorname{supp}(x)$. Therefore
\[
I=\bigcup_{x\in E}\operatorname{supp}(x).
\]
Every support in this union is finite. Since $I$ is infinite, it follows that $|I|\leq|E|$. Thus
\[
d\left(\bigoplus_{i\in I}K_i\right)=|I|.\qedhere
\]
\end{proof}


\section{Proof of Theorem \ref{ThA}}\label{s2}

In this section, we prove Theorem \ref{ThA}. The finite case is immediate: a finite compact group is discrete, so its only dense subgroup is the whole group and $d(G)=w(G)=|G|$. Hence we may restrict attention to infinite compact groups. We first establish property $\JM^*$ for all infinite profinite groups. We then turn to arbitrary infinite compact groups and use the structure of the identity component to reduce the remaining case to invariant dense subgroups of connected abelian and semisimple compact groups.

\subsection{The Profinite Case}

We begin with the totally disconnected case. Recall that a subset $X$ of a profinite group $G$ is called a \emph{generating set converging to the identity} if $\overline{\langle X\rangle}=G$
and, for every open neighbourhood $U$ of the identity $e$, the set $X\setminus U$ is finite. Such sets are also called \emph{suitable generating sets}. If $X$ is infinite, this is equivalent to saying that $X\setminus\{e\}$ is discrete and that, after adjoining $e$ if necessary, $X\cup\{e\}$ is the one-point compactification of $X\setminus\{e\}$, with $e$ as its unique accumulation point. If $X$ is finite, the convergence condition is automatic, and the existence of such an $X$ simply means that $G$ is topologically finitely generated.

Every profinite group admits a suitable generating set. We denote by $\delta(G)$ the least cardinality of such a set. One has
\begin{lemma}\label{delta=w}\cite[Proposition 2.6.1 and Corollary 2.6.3]{RZ}
If $G$ is a profinite group and $\delta(G)$ is infinite, then $\delta(G)=w(G)$.
\end{lemma}

For later use, we also fix the standard meaning of rank for free pro-$p$ groups; in this paper, the notation $\operatorname{rank}(F)$ will be used only when $F$ is free pro-$p$. This convention should not be confused with other notions of rank for general profinite or pro-$p$ groups. If $F=F_p(X)$ is a free pro-$p$ group on a set $X$ converging to the identity, then
\[
\operatorname{rank}(F):=|X|.
\]
This is well defined by \cite[Lemma~3.3.5]{RZ}. Moreover, with the notation introduced above, every free pro-$p$ group satisfies
\[
\operatorname{rank}(F)=\delta(F).
\]
Indeed, if $X$ is finite, then $X$ is a suitable generating set, while the canonical Frattini quotient
\[
F/\Phi(F)\cong(\mathbb Z/p\mathbb Z)^X
\]
shows that no set of fewer than $|X|$ elements can topologically generate $F$. Hence $\delta(F)=|X|$. If $X$ is infinite, the same quotient shows that $F$ is not topologically finitely generated, so $\delta(F)$ is infinite. By \cite[Proposition~2.6.2(b)]{RZ}, $w(F)=|X|$, while Lemma~\ref{delta=w} gives $\delta(F)=w(F)$. Thus again $\delta(F)=|X|=\operatorname{rank}(F)$.

The following observation settles the upper endpoint of the density spectrum for arbitrary profinite groups.

\begin{proposition}\label{maxdense}
Let $G$ be an infinite profinite group. Then there exists a dense subgroup $H\leq G$ such that $d(H)=w(G)$.
\end{proposition}

\begin{proof}
Put $\kappa=w(G)$. First suppose that $\kappa=\omega$. Then $G$ is compact metrizable, so it has a countable dense subset $D$. The subgroup $H=\langle D\rangle$ is countable and dense in $G$. Since $H$ is infinite and Hausdorff, it has no finite dense subset. Hence $d(H)=\omega=\kappa$.

Now suppose that $\kappa>\omega$. We first observe that $\delta(G)$ is infinite. Indeed, if $\delta(G)$ were finite, then $G$ would be metrizable, contradicting $\kappa>\omega$.
Therefore Lemma~\ref{delta=w} gives $\delta(G)=w(G)=\kappa$. Choose a suitable generating set $X\subseteq G$ with $|X|=\kappa$, and put $H=\langle X\rangle$. Then $H$ is dense in $G$, and clearly
$d(H)\leq |H|\leq\kappa$.

Suppose, towards a contradiction, that $d(H)=\lambda<\kappa$. Choose a dense set $D\subseteq H$ with $|D|=\lambda$. For every $d\in D$, choose a finite set $F_d\subseteq X$ such that $d\in\langle F_d\rangle$, and put $Y=\bigcup_{d\in D}F_d$.
Then $D\subseteq\langle Y\rangle$ and $|Y|\leq\max\{\lambda,\omega\}<\kappa$. Since $D$ is dense in $H$ and $H$ is dense in $G$, we have $\overline{D}=G$. Therefore
$G=\overline{D}\subseteq\overline{\langle Y\rangle}$, and hence
$\overline{\langle Y\rangle}=G$.

Finally, $Y\subseteq X$, so $Y$ still converges to $e$: for every neighbourhood $U$ of $e$,
$Y\setminus U\subseteq X\setminus U$, and the latter set is finite. Thus $Y$ is a suitable generating set for $G$ of cardinality strictly smaller than $\kappa=\delta(G)$, a contradiction. Therefore $d(H)\geq\kappa$, and consequently $d(H)=\kappa$.
\end{proof}

We next recall the quotient lemma needed to obtain all intermediate densities.

\begin{lemma}\label{weight_quotient}
Let $G$ be an infinite profinite group. Then for every infinite cardinal $\kappa\leq w(G)$, there exists a closed normal subgroup $N$ of $G$ such that $w(G/N)=\kappa$.
\end{lemma}

\begin{proof}
Put $\tau=w(G)$. Since the open normal subgroups form a local base at the identity of a profinite group, there is a family of $\tau$ distinct open normal subgroups of $G$. Choose distinct open normal subgroups $U_\alpha$, $\alpha<\kappa$, and set
\[
N:=\bigcap_{\alpha<\kappa}U_\alpha.
\]
Then $N$ is a closed normal subgroup of $G$. The diagonal homomorphism
\[
G/N\longrightarrow \prod_{\alpha<\kappa}G/U_\alpha
\]
is a topological embedding. Since each factor $G/U_\alpha$ is finite, it follows that $w(G/N)\leq\kappa$.

For the reverse inequality, the subgroups $U_\alpha/N$, $\alpha<\kappa$, are $\kappa$ distinct open normal subgroups of $G/N$. Suppose that $w(G/N)=\mu<\kappa$. The quotient $G/N$ is infinite, so $\mu$ is infinite. As $G/N$ is profinite, it has a local base $\mathcal B$ at the identity consisting of open normal subgroups with $|\mathcal B|=\mu$. Every open normal subgroup $V$ of $G/N$ contains some $B\in\mathcal B$. For a fixed $B\in\mathcal B$, the open normal subgroups containing $B$ correspond to normal subgroups of the finite group $(G/N)/B$, and hence there are only finitely many of them. Therefore $G/N$ has at most $\mu$ open normal subgroups, contradicting the existence of the distinct family $\{U_\alpha/N:\alpha<\kappa\}$. Thus $w(G/N)\geq\kappa$, and consequently $w(G/N)=\kappa$.
\end{proof}

\begin{proposition}\label{profiniteJM}
Every infinite profinite group satisfies property $\JM^*$.
\end{proposition}

\begin{proof}
Let $G$ be an infinite profinite group, and fix a cardinal
\[
\kappa\in[d(G),w(G)].
\]
By Lemma~\ref{weight_quotient}, there exists a closed normal subgroup $N$ of $G$ such that the quotient $Q=G/N$ has weight $w(Q)=\kappa$. Since $\kappa\geq d(G)$, we also have $d(Q)\leq d(G)\leq\kappa$.

By Proposition~\ref{maxdense}, the profinite group $Q$ has a dense subgroup $H_0$ with
\[
d(H_0)=w(Q)=\kappa.
\]
Let $\pi:G\to Q$ be the canonical projection. Lemma~\ref{preimage} now shows that $H=\pi^{-1}(H_0)$ is a dense subgroup of $G$ satisfying $d(H)=\kappa$. Since $\kappa$ was arbitrary, $G$ satisfies property $\JM^*$.
\end{proof}

\subsection{The Group of Automorphisms}

A natural object associated with any topological group $G$ is its \emph{automorphism group} $\Aut(G)$, consisting of all topological group automorphisms of $G$. When $G$ is locally compact, we endow $\Aut(G)$ with the \emph{Birkhoff topology}. This topology is generated by a subbasis consisting of sets of the form
\[
\mathcal{B}(L, U) := \{ \phi \in \Aut(G) : \phi(x) \in Ux \text{ and } \phi^{-1}(x) \in Ux \text{ for all } x \in L \},
\]
where $L$ is a compact subset of $G$, and $U$ is an open neighbourhood of the identity in $G$. Originally introduced by Birkhoff in the context of uniform structures, this topology is specifically designed to reflect both the topological and algebraic features of $\Aut(G)$.

Under the Birkhoff topology, $\Aut(G)$ becomes a topological group. Furthermore, if a topological group $K$ acts on $G$ via topological automorphisms, then this action is continuous if and only if the induced group homomorphism $K \to \Aut(G)$ is continuous. For further details, we refer the reader to \cite[p.~260]{DPS}.

\subsection{Dense Subgroups Invariant under a Compact Group Action}
\subsubsection{The Abelian Case}

In this part, all group operations are written additively, with $0$ denoting the identity element.

Let $G$ be a locally compact abelian group, and let $X$ denote its Pontryagin dual group. It is a classical fact that $\Aut(G)$ and $\Aut(X)$ are \textbf{topologically anti-isomorphic} (see \cite[(26.9) Theorem]{HR}). More precisely, the map
\[
\ad: \Aut(G) \to \Aut(X), \quad \sigma \mapsto \widehat{\sigma}
\]
is both a group anti-isomorphism and a homeomorphism.

The following lemma is essentially folklore, though elementary; we include its proof for the sake of completeness.

\begin{lemma}\label{dual}
Let $G$ be a locally compact abelian group and let $\sigma \in \Aut(G)$. Then for any closed subgroup $H \leq G$, the subgroup $H$ is $\sigma$-invariant if and only if its annihilator $H^\perp$ is $\widehat{\sigma}$-invariant.
\end{lemma}

\begin{proof}
Suppose $H \leq G$ is $\sigma$-invariant. Then for any $\chi \in H^\perp$ and $h \in H$, we have
\[
\widehat{\sigma}(\chi)(h) = \chi(\sigma(h)) \in \chi(H) = \{0\},
\]
which implies $\widehat{\sigma}(\chi) \in H^\perp$. Thus, $H^\perp$ is $\widehat{\sigma}$-invariant.

The converse follows symmetrically by Pontryagin duality.
\end{proof}

\begin{remark}\label{action}
By convention, when we speak of an ``action'' without qualification, we mean a left action. The canonical anti-isomorphism $\sigma \mapsto \widehat{\sigma}$ naturally yields a \emph{right} action of $\Aut(G)$ on the dual group $X$, given by $\chi \cdot \sigma := \widehat{\sigma}(\chi) = \chi \circ \sigma$. Consequently, any continuous left action of a topological group $K$ on a locally compact abelian group $G$ inherently induces a continuous right action on its dual $X$.

To maintain notational consistency with left modules in subsequent proofs, we canonically convert this right action into a left action by composing with the inverse. That is, the induced left action of $K$ on $X$ is defined by
\[
k \cdot \chi := \chi \circ k^{-1} = \widehat{k^{-1}}(\chi).
\]
In what follows, we will freely view the dual group $X$ as being equipped with this induced continuous left action whenever $G$ admits a continuous left $K$-action. Because a closed subgroup of $G$ is $\sigma$-invariant if and only if it is $\sigma^{-1}$-invariant, Lemma \ref{dual} can be elegantly restated in terms of this induced left action: \emph{a closed subgroup $H \leq G$ is $\sigma$-invariant if and only if its annihilator $H^\perp$ is $\sigma$-invariant.}
\end{remark}

\begin{lemma}\label{abinv}
Let $G$ be an infinite connected compact abelian group, and let $K$ be a compact subgroup of $\Aut(G)$. Then for every cardinal $\kappa \in [d(G), w(G)]$, there exists a dense $K$-invariant subgroup $H \leq G$ such that $d(H) = \kappa$.
\end{lemma}

\begin{proof}
Without loss of generality, we may assume that $w(G)$ is uncountable. We can also consider only the case that $\kappa>\omega$ because otherwise $\kappa=d(G)=\omega$ and we trivially take $H=G$.  Let $A := \widehat{G}$ be the dual group of $G$. Since $G$ is compact and connected, $A$ is a torsion-free discrete abelian group. By the preceding discussion, the left action of $K$ on $G$ canonically induces a continuous left action of $K$ on $A$. 

Let $V=A\otimes_{\mathbb{Z}}\mathbb{Q}$ be the divisible hull of $A$ (see \cite[Sec. 24]{Fuc}). Then $V$ is a $\mathbb{Q}$-vector space. It is a standard fact that every automorphism of $A$ uniquely extends to a $\mathbb{Q}$-linear automorphism of $V$. Hence, the continuous action of $K$ on $A$ naturally extends to a $\mathbb{Q}$-linear representation of $K$ on $V$, such that $A$ becomes a $K$-invariant subgroup of $V$.

Because $A$ is discrete and the action of $K$ is continuous, the stabilizer of any element in $A$ is an open subgroup of the compact group $K$. Consequently, the $K$-orbit $K \cdot x$ is finite for every $x \in A$.

We shall construct, by transfinite induction, a family $\{A_\alpha : \alpha < \kappa\}$ of non-zero $K$-invariant countable subgroups of $A$ such that
\[ A_\alpha \cap \sum_{\beta < \alpha} A_\beta = \{0\} \quad \text{for all } \alpha < \kappa. \]

We begin by choosing a non-zero element $a_0 \in A$. Let $V_0$ be the $\mathbb{Q}$-linear subspace of $V$ generated by the finite orbit $K \cdot a_0$. Then $V_0$ is $K$-invariant and finite-dimensional, hence countable. Define $A_0 := A \cap V_0$, which is a non-zero $K$-invariant subgroup of $A$.

Now suppose that for some ordinal $\alpha < \kappa$, we have successfully constructed the family $\{A_\beta\}_{\beta < \alpha}$. Let $V_\alpha$ be the $\mathbb{Q}$-linear span of $\sum_{\beta < \alpha} A_\beta$ in $V$. Since $V_\alpha$ has cardinality at most $\omega \cdot |\alpha| < \kappa\leq w(G) = |A|$, we can pick an element $a \in A \setminus V_\alpha$.

Let $W$ be the $\mathbb{Q}$-subspace of $V$ generated by the finite orbit $K \cdot a$. Then $W$ is $K$-invariant and finite-dimensional. Since $V_\alpha$ is also $K$-invariant, the intersection $U := W \cap V_\alpha$ is a $K$-invariant subspace of $W$, and $\dim U < \dim W$ (because $a \notin U$).

Let $M := \{k \in K : k \cdot x = x \text{ for all } x \in K \cdot a\}$ be the pointwise stabilizer of the finite set $K \cdot a$. Since $K$ acts continuously on the discrete space $A$, $M$ is an open normal subgroup of $K$. Thus, the quotient $K/M$ is a finite group. Since $M$ acts trivially on $W$, the action of $K$ on $W$ factors through the finite group $K/M$. 

By Maschke's Theorem (see e.g., \cite[p. 116]{AB}), since the underlying field $\mathbb{Q}$ has characteristic zero and $K/M$ is finite, the $K/M$-invariant subspace $U$ admits a $K/M$-invariant (and thus $K$-invariant) complement $W_\alpha \leq W$ such that $W = U \oplus W_\alpha$. 

Define $A_\alpha := A \cap W_\alpha$. Because $\dim U < \dim W$, we have $\dim W_\alpha > 0$. Since $V=A\otimes_{\mathbb Z}\mathbb Q$, every non-zero vector $v\in W_\alpha$ has a non-zero integer multiple lying in $A$. Hence $W_\alpha\cap A\neq\{0\}$, and therefore $A_\alpha$ is a non-zero $K$-invariant countable subgroup of $A$. This completes the induction.

Let $S := \bigoplus_{\alpha < \kappa} A_\alpha$ denote the resulting internal direct sum in $A$. Since $S$ is $K$-invariant, Lemma \ref{dual} ensures that its annihilator $N := S^\perp \leq G$ is $K$-invariant. Furthermore, we have the natural identification
\[ G/N \cong \widehat{S} \cong \prod_{\alpha < \kappa} \widehat{A_\alpha}. \]
Let $\pi: G \to \prod_{\alpha < \kappa} \widehat{A_\alpha}$ denote the canonical projection. Since $N$ is $K$-invariant, this projection is $K$-equivariant. For each $\alpha < \kappa$, defining $S_\alpha := \bigoplus_{\beta \neq \alpha} A_\beta$ and $N_\alpha := S_\alpha^\perp$ shows that $N_\alpha = \pi^{-1}(\widehat{A_\alpha})$, where $\widehat{A_\alpha}$ is identified with the corresponding coordinate subgroup of the product. Since each $S_\alpha$ is $K$-invariant, each $N_\alpha$ is $K$-invariant; by equivariance of $\pi$, the coordinate subgroup $\widehat{A_\alpha}$ is therefore $K$-invariant as well.

Let $P$ be the finite-support subgroup $\bigoplus_{\alpha<\kappa}\widehat{A_\alpha}$ inside $\prod_{\alpha < \kappa} \widehat{A_\alpha}$. Since each $\widehat{A_\alpha}$ is a non-trivial compact metrizable group, Lemma \ref{prod} guarantees that $P$ is a dense subgroup of $G/N$ with density $\kappa$. Moreover, since the $K$-action preserves the coordinates, $P$ is $K$-invariant.

Finally, define $H := \pi^{-1}(P)$. Then $H$ is a dense $K$-invariant subgroup of $G$. Since $d(N) \leq d(G) \leq \kappa$, it follows from Lemma \ref{preimage} that $d(H) = \kappa$. The proof is complete.
\end{proof}

\subsubsection{Semisimple Connected Compact Groups}

By a \emph{semisimple connected compact group}, we mean a connected compact group whose commutator subgroup coincides with the whole group. For a detailed exposition on their structure, we refer to \cite[Ch.~9]{HM}.

\begin{lemma}\label{semisim}
Let $G$ be a semisimple connected compact group, and let $K$ be a compact subgroup of $\Aut(G)$. Then for any cardinal $\kappa \in [d(G), w(G)]$, there exists a dense $K$-invariant subgroup $H \leq G$ such that $d(H) = \kappa$.
\end{lemma}

\begin{proof}
Without loss of generality, we may assume that $G$ is not metrizable.

By the structure theory of semisimple connected compact groups \cite[Theorem~9.19]{HM}, the centre $Z$ of $G$ is totally disconnected, and there exists a topological isomorphism
\[
G/Z \cong \prod_{\alpha < \tau} S_\alpha,
\]
where each $S_\alpha$ is a simple, connected, compact Lie group with trivial centre, and $\tau = w(G)$. 
It follows that $w(G/Z) = w(G)$ and $d(G/Z) = d(G)$.

Moreover, since the centre $Z$ is strictly invariant under any automorphism, it is clearly $K$-invariant. Hence, the action of $K$ on $G$ descends to a continuous action on $G/Z$. 
Therefore, it suffices to construct a $K$-invariant dense subgroup $H_0 \leq G/Z$ of density $\kappa$, since its preimage $H := \pi^{-1}(H_0)$ under the canonical projection $\pi: G \to G/Z$ will then be a dense $K$-invariant subgroup of $G$ with $d(H) = \kappa$, by Lemma~\ref{preimage}.

Thus, we may reduce our problem to the case where $G = \prod_{i \in I} S_i$, with each $S_i$ being a simple, connected, compact Lie group with trivial centre, and $|I| = \tau = w(G)$.

It is a standard fact that the connected normal subgroups of $G$ are precisely the subproducts over subsets of $I$. 
Hence, for each $\sigma \in K$ and $i \in I$, there exists a unique $j \in I$ such that $\sigma(S_i) = S_j$. This yields a well-defined homomorphism
\[
\varphi: K \to \Sym(I),
\]
where $\Sym(I)$ denotes the group of all permutations of $I$, endowed with the topology of pointwise convergence (treating $I$ as a discrete space).

We claim that $\varphi$ is continuous. 

Since a subbasis at the identity of $\Sym(I)$ consists of point stabilizers, it suffices to show that for each $i \in I$, the setwise stabilizer of $S_i$ in $K$ is an open subgroup.

Let $L := \prod_{j \in I \setminus \{i\}} S_j$. Being a closed subgroup of the compact group $G$, $L$ is compact. Since $S_i$ is a Lie group, it possesses the ``no small subgroups'' (NSS) property; thus, we can choose an open neighbourhood $V$ of the identity in $S_i$ that contains no non-trivial subgroups. 
Then $U := L \times V$ is an open neighbourhood of the identity in $G$. 

Now consider the basic open set in the Birkhoff topology of $\Aut(G)$:
\[
\mathcal{B}(L, U) := \{\sigma \in K : \sigma(x) \in Ux \text{ and } \sigma^{-1}(x) \in Ux \text{ for all } x \in L\}.
\]
We claim that each $\sigma \in \mathcal{B}(L, U)$ satisfies $\sigma(L) = L$, and consequently $\sigma(S_i) = S_i$.

Indeed, for any $x \in L$, we have $\sigma(x) \in Ux \subseteq UL = U$, which implies $\sigma(L) \subseteq U$. 
However, all subgroups of $G$ contained entirely within $U = L \times V$ are contained in $L$ (as $V$ contains no non-trivial subgroups). Therefore, we must have $\sigma(L) \leq L$. 
Applying the identical argument to $\sigma^{-1}$ yields $\sigma^{-1}(L) \leq L$. Together, these force $\sigma(L)=L$. Since all factors have trivial centre, $C_G(L)=S_i$. Hence
$\sigma(S_i)=\sigma(C_G(L))=C_G(\sigma(L))=C_G(L)=S_i$.

This proves that the stabilizer of $i$ in $K$ contains the open neighbourhood $\mathcal{B}(L, U)$, establishing the continuity of $\varphi$.

Because $K$ is compact and $\varphi$ is continuous, $\varphi(K)$ is a compact subgroup of $\Sym(I)$. 
Since the action of $\varphi(K)$ on the discrete space $I$ is continuous, the stabilizer of each $i \in I$ is an open subgroup of $\varphi(K)$. Compactness then dictates that this stabilizer has finite index, which guarantees that the orbit $\varphi(K) \cdot i$ is finite. This yields a partition of $I$ into finite $K$-orbits:
\[
I = \bigcup_{\alpha < \tau} I_\alpha,
\]
where each $I_\alpha$ is a finite $K$-orbit.

For each $\alpha<\tau$, set $T_\alpha=\prod_{i\in I_\alpha}S_i$. Therefore $G=\prod_{\alpha<\tau}T_\alpha$.

Now fix any cardinal $\kappa \in [d(G), \tau]$. Set
$A=\prod_{\kappa\leq\alpha<\tau}T_\alpha$, and let $B$ be the finite-support subgroup $\bigoplus_{\beta<\kappa}T_\beta$ inside $\prod_{\beta<\kappa}T_\beta$.

Because each orbit $I_\alpha$ is finite, the block subgroup $T_\alpha = \prod_{i \in I_\alpha} S_i$ is a finite product of compact metrizable Lie groups, hence itself compact and metrizable. By Lemma~\ref{prod}, the finite-support subgroup $B$ is a dense subgroup of the corresponding product space and has density precisely $\kappa$.

The projection of $G$ onto $A$ is continuous and surjective, so $d(A)\leq d(G)\leq\kappa$. Hence
\[
\kappa=d(B)\leq d(A\times B)\leq d(A)\cdot d(B)\leq\kappa,
\]
and therefore $d(A\times B)=\kappa$.
Because both $A$ and $B$ are constructed from unions of complete $K$-orbits $I_\alpha$, they are strictly $K$-invariant. Thus, $H_0 := A \times B$ is a dense $K$-invariant subgroup of $G$ with density $\kappa$, concluding the proof.
\end{proof}

\subsection{The General Compact Case}

We are now in a position to complete the proof of Theorem \ref{ThA}. 

The argument relies heavily on the structure theory of compact groups, the details of which can be found in \cite[Ch.~9]{HM}. 
Throughout this subsection, let $G$ be an infinite compact group, and let $G_0$ denote the connected component of the identity in $G$.

By standard cardinal arithmetic for topological groups, we have the identity $w(G) = w(G/G_0) \cdot w(G_0)$. Since $w(G)$ is infinite, it immediately follows that either $w(G) = w(G/G_0)$ or $w(G) = w(G_0)$. We analyze these two cases separately.

\begin{case}\label{case1}
$w(G/G_0) = w(G)$.
\end{case}

\begin{proof}
In this case, we naturally have $d(G/G_0) = d(G) = \log w(G)$. 
 Since $G/G_0$ is profinite, the preceding proposition shows that $G/G_0$ satisfies $\JM^*$.

For any infinite cardinal
\[
\kappa \in [d(G), w(G)] = [d(G/G_0), w(G/G_0)],
\]
the property $\JM^*$ yields a dense subgroup $H \leq G/G_0$ with density $d(H) = \kappa$. 
By Lemma~\ref{preimage}, the preimage of $H$ under the canonical continuous open projection $\pi: G \to G/G_0$ forms a dense subgroup of $G$ with density exactly $\kappa$. This completes the proof for Case \ref{case1}.
\end{proof}

\begin{case}\label{case2}
$w(G_0) = w(G)$.
\end{case}

\begin{proof}
We begin by systematically reducing the ambient group. By the Compact-group Sandwich Theorem, in the form \cite[Corollary~9.42]{HM}, there exists a closed totally disconnected subgroup $\widetilde D\leq G$ such that $G=G_0\widetilde D$. Put
\[
\Delta=G_0\cap\widetilde D.
\]
Then $\Delta$ is a closed totally disconnected abelian normal subgroup of $G$, and
\[
G/\Delta\cong (G_0/\Delta)\rtimes (\widetilde D/\Delta).
\]
In particular, $E:=\widetilde D/\Delta$ is profinite. Notice that centrality of $\Delta$ is needed only inside the connected group $G_0$: since $\Delta$ is totally disconnected and normal in $G_0$, it is central in $G_0$.

By Lemma~\ref{conq},
\[
w(G_0/\Delta)=w(G_0)=w(G).
\]
Since $G_0/\Delta$ is a closed subgroup of $G/\Delta$, we obtain
\[
w(G)\geq w(G/\Delta)\geq w(G_0/\Delta)=w(G),
\]
and hence $w(G/\Delta)=w(G)$. Both groups are infinite compact groups, so $d(G/\Delta)=d(G)$. By Lemma~\ref{preimage}, it therefore suffices to prove property $\JM^*$ for $G/\Delta$. Thus, after replacing $G$ by this quotient, we may assume that
\[
G=G_0\rtimes E
\]
with $E$ profinite.

Let $G_0'$ denote the commutator subgroup of $G_0$, which is closed and semisimple. 
By \cite[Theorem~9.24]{HM}, we can write $G_0 = A G_0'$, where $A$ denotes the connected component of the identity in the centre $Z(G_0)$. Their intersection $\Gamma := A \cap G_0'$ is totally disconnected, and it follows that
\[
G_0/\Gamma \cong (A/\Gamma) \times (G_0'/\Gamma).
\]
Since $A$ and $G_0'$ are topologically characteristic subgroups of $G_0$, and $G_0$ is characteristic in $G$, it follows that $\Gamma$ is normal in $G$. Applying Lemma~\ref{conq} once more, we obtain
\[
w(G)\geq w(G/\Gamma)\geq w(G_0/\Gamma)=w(G_0)=w(G).
\]
Thus $w(G/\Gamma)=w(G)$ and, since the groups are infinite and compact, $d(G/\Gamma)=d(G)$. By Lemma~\ref{preimage} again, it suffices to prove the result for the quotient $G/\Gamma$. 
Thus, we may further reduce the problem by assuming that $G_0 \cong A \times S$, where $A$ is an abelian connected compact group and $S$ is a semisimple connected compact group, such that
\[
G \cong (A \times S) \rtimes E,
\]
where $A$ and $S$ are both normal in $G$, and the profinite group $E$ acts continuously on each.

Because $G_0 \cong A \times S$, we have either $w(A) = w(G_0) = w(G)$ or $w(S) = w(G_0) = w(G)$. 

Assume first that $w(A) = w(G)$. 
Fix an arbitrary cardinal $\kappa \in [d(G), w(G)]$. 
By Lemma~\ref{abinv}, there exists a dense $E$-invariant subgroup $H \leq A$ with $d(H) = \kappa$. 
Because $H$ is $E$-invariant, the direct product $H \times S$ is naturally $E$-invariant, allowing us to form the semidirect product
\[
H_1 := (H \times S) \rtimes E.
\]
Because $H$ is dense in $A$, $H_1$ is a dense subgroup of $G$. As a space, $H_1$ is homeomorphic to $H\times S\times E$. Since the projection onto $H$ is continuous and surjective, while $d(S),d(E)\leq d(G)\leq\kappa$, we obtain
\[
\kappa=d(H)\leq d(H_1)\leq d(H)\cdot d(S)\cdot d(E)\leq\kappa.
\]
Thus $d(H_1)=\kappa$.

If, instead, $w(S) = w(G)$, we apply Lemma~\ref{semisim} to construct a dense $E$-invariant subgroup of $S$ with density $\kappa$. An identical semidirect product construction then yields the desired dense subgroup of $G$.
\end{proof}

Combining Cases~\ref{case1} and \ref{case2} establishes property $\JM^*$ for $G$, completing the proof of Theorem \ref{ThA}.

\section{Density of Closed Subgroups}\label{s3}

In this section, we turn our attention to another important class of subgroups: closed subgroups. We begin our discussion with separable topological groups. It was shown by Comfort and Itzkowitz \cite{CIt} that a closed subgroup of a locally compact group $G$ has density at most $d(G)$. Consequently, closed subgroups of separable locally compact groups---and, in particular, of separable compact groups---are themselves separable.

To understand the behavior of closed subgroups beyond the locally compact realm, it is helpful to briefly recall a hierarchy of compactness-like properties. A Tychonoff space $X$ is called \emph{pseudocompact} if every continuous real-valued function on $X$ is bounded, and \emph{countably compact} if every countable open cover of $X$ admits a finite subcover. It is a standard fact that every countably compact space is pseudocompact. In the context of topological groups, a group is called \emph{precompact} (or totally bounded) if it is topologically isomorphic to a subgroup of a compact group. While every pseudocompact group is precompact, the converse fails in general.

Closed subgroups behave quite differently beyond the locally compact setting. Leiderman, Morris, and Tkachenko \cite{LMT} established that the density of a closed subgroup of a separable precompact, or even pseudocompact, group can be arbitrarily large. More precisely, they proved the following result:

\begin{theorem}[Leiderman--Morris--Tkachenko]\label{thm_LMT}
Every precompact group of weight at most $\mathfrak{c}$ is topologically isomorphic to a closed subgroup of a separable, connected, pseudocompact group of weight at most $\mathfrak{c}$. In the abelian case, the ambient group can also be chosen abelian.
\end{theorem}

Recall that subspaces of separable Tychonoff spaces have weight at most $\mathfrak{c}$, which follows from a well-known cardinal inequality. Thus, even a separable pseudocompact group may contain a closed subgroup of density $\mathfrak{c}$, the largest value allowed by this inequality.

Leiderman, Morris, and Tkachenko also investigated countably compact groups in the same paper. A key distinction between countable compactness and pseudocompactness is that countable compactness is inherited by closed subspaces. Therefore, their general embedding result cannot extend to countably compact groups, as not every precompact group is countably compact.

Nevertheless, assuming the Continuum Hypothesis ($\mathbf{CH}$), they constructed a separable countably compact group containing a closed subgroup of density $\mathfrak{c}$. Their construction uses an $\omega$-HFD, whose existence follows from $\mathbf{CH}$. They concluded with the following open question:

\begin{question}\cite[Question 5.6]{LMT}
Does there exist in $\mathbf{ZFC}$ a separable countably compact group with a closed non-separable subgroup?
\end{question}

The main objective of this section is to answer this question affirmatively in $\mathbf{ZFC}$.

\subsection{Proof of Theorem \ref{ThB}}
Throughout this subsection, all groups are abelian. We write the group operation additively and denote the identity element by $0$.

We first recall a recent result on countably compact groups. A long-standing problem asked whether $\mathbf{ZFC}$ proves the existence of a countably compact group without non-trivial convergent sequences.
An affirmative answer to this question would resolve several conjectures concerning countably compact groups, such as the existence of a countably compact group whose square is not countably compact. The question was answered by Hru\v{s}\'{a}k, van Mill, Ramos-Garc\'{i}a, and Shelah \cite{HMRS}, who provided a $\mathbf{ZFC}$ construction of such a group.

Before stating their breakthrough result, we briefly recall the concept of the Bohr topology. For a discrete abelian group $A$, the \emph{Bohr topology} on $A$, often denoted by $A^\#$, is the maximal precompact (totally bounded) group topology on $A$. Equivalently, it is the initial topology induced by the family of all homomorphisms from $A$ into the circle group $\mathbb{T}$. Every homomorphism from a Boolean group into $\mathbb{T}$ has image contained in the unique subgroup of order $2$, which we identify with $\mathbb{Z}(2)$.  We write $2=\{0,1\}\cong\mathbb{Z}(2)$. Consequently, the Bohr topology of a Boolean group is generated by all homomorphisms into $2$.

Their construction begins with a countable Boolean group equipped with the Bohr topology, and subsequently extends both its algebraic structure and its topology to achieve the desired properties. This yields the following crucial lemma:

\begin{lemma}\label{Le0}
In $\mathbf{ZFC}$, there exists a countably compact Boolean group $G$ of cardinality $\mathfrak{c}$ with a countably infinite dense subgroup $G_0$, such that the induced topology on $G_0$ is the Bohr topology.
\end{lemma}

\begin{proof}
We extract the required additional information from the construction in the proof of \cite[Theorem~4.1]{HMRS}. Let
$A=[\omega]^{<\omega}$ and $B=[\mathfrak c]^{<\omega}$, viewed as Boolean groups under symmetric difference. In that construction, every homomorphism $\Phi:A\to 2$ is extended to a homomorphism $\overline\Phi:B\to 2$, and $B$ is endowed with the group topology generated by the family of all these extensions. If
\[
K=\bigcap_{\Phi\in\operatorname{Hom}(A,2)}\ker(\overline\Phi),
\]
then the proof of \cite[Theorem~4.1]{HMRS} shows that $K$ is finite and that the Hausdorff quotient
$G=B/K$ is countably compact. Hence $|G|=\mathfrak c$.

Because $\overline\Phi|_A=\Phi$ for every $\Phi\in\operatorname{Hom}(A,2)$, we have $K\cap A=\{0\}$. Thus the image $G_0$ of $A$ in $G$ is a countably infinite subgroup, and the topology induced on $G_0$ is exactly the Bohr topology of $A$. Finally, consider the evaluation embedding
\[
e:G\longrightarrow 2^{\operatorname{Hom}(A,2)},
\qquad
e(b+K)(\Phi)=\overline\Phi(b).
\]
We show directly that $e(G_0)$ is dense in $e(G)$. Fix $b\in B$ and finitely many homomorphisms $\Phi_1,\ldots,\Phi_n\in\operatorname{Hom}(A,2)$. It is enough to find $a\in A$ such that
\[
\Phi_j(a)=\overline\Phi_j(b)\qquad (j=1,\ldots,n).
\]
A transfinite induction through the recursive extension procedure in \cite[Theorem~4.1]{HMRS}, using the additivity of ultrafilter limits in the finite discrete group $2$ at each extension step, shows that the assignment $\Phi\mapsto\overline\Phi$ is $\mathbb F_2$-linear. In particular, if $\Phi,\Psi\in\operatorname{Hom}(A,2)$, then
$\overline{\Phi+\Psi}=\overline\Phi+\overline\Psi$. Hence every linear relation among $\Phi_1,\ldots,\Phi_n$ is satisfied by the tuple
$(\overline\Phi_1(b),\ldots,\overline\Phi_n(b))$. By elementary linear algebra over $\mathbb F_2$, this tuple therefore lies in the image of the evaluation map
\[
A\longrightarrow 2^n,
\qquad
a\longmapsto(\Phi_1(a),\ldots,\Phi_n(a)).
\]
Thus such an $a$ exists. Every basic neighbourhood of $e(b+K)$ in $e(G)$ consequently meets $e(G_0)$, proving that $G_0$ is dense in $G$.
\end{proof}

Let $G$ and $G_0$ be as in Lemma~\ref{Le0}. Let $X$ denote the group of continuous characters of $G$. Since $G_0$ is dense in $G$, the two groups have the same Raikov completion. Every continuous character of $G_0$ extends uniquely to this completion, and hence to $G$, while uniqueness on the dense subgroup gives injectivity of restriction. Therefore restriction induces an isomorphism between the continuous character groups of $G$ and $G_0$. Since $G_0$ is a countably infinite Boolean group endowed with its Bohr topology, we have $X\cong 2^\omega$, and hence $\dim_{\mathbb F_2}X=|X|=\mathfrak c$.

For each nonzero element $g\in G$, let
\[
X_g=\{\chi\in X:\chi(g)=1\}.
\]

The following lemma is straightforward.

\begin{lemma}\label{Le1}
For every subgroup $Y\le X$ with $|Y|<\mathfrak c$ and every nonzero element $g\in G$, there exist linearly independent over $\mathbb{F}_2$ elements $e^{(1)},e^{(2)}\in X_g$ such that
\[
Y\cap\langle e^{(1)},e^{(2)}\rangle=\{0\}.
\]
\end{lemma}

\begin{proof}
Since $G$ is Hausdorff and precompact, its continuous characters separate points. Hence $X_g\neq\emptyset$. Moreover, $X_g$ generates $X$: indeed, if $\chi_0\in X_g$ and $\chi\in X\setminus X_g$, then $\chi+\chi_0\in X_g$, and hence $\chi=\chi_0+(\chi+\chi_0)$. Since $\dim_{\mathbb F_2}X=\mathfrak c$ and $|Y|<\mathfrak c$, the
image of $X_g$ in $X/Y$ spans $X/Y$, which has dimension
$\mathfrak c$. Hence we may choose $e^{(1)},e^{(2)}\in X_g$ whose images in $X/Y$ are linearly independent. This gives $Y\cap\langle e^{(1)},e^{(2)}\rangle=\{0\}$.
\end{proof}

\begin{proof}[Proof of Theorem \ref{ThB}] Fix an enumeration
\[
G\setminus\{0\}=\{g_\alpha:\alpha<\mathfrak c\}.
\]

We shall construct, by transfinite induction, elements
\[
e_\alpha^{(1)},e_\alpha^{(2)}\in X_{g_\alpha}
\qquad (\alpha<\mathfrak c)
\]
such that the set
\[
\{e_\alpha^{(i)}:\alpha<\mathfrak c,\ i=1,2\}
\]
is linearly independent over $\mathbb F_2$.

To begin, apply Lemma~\ref{Le1} with $Y=\{0\}$ and $g=g_0$ to choose linearly independent elements $e_0^{(1)},e_0^{(2)}\in X_{g_0}$. Assume that for some $0<\alpha<\mathfrak c$, we have already chosen
\[
\{e_\beta^{(1)},e_\beta^{(2)}:\beta<\alpha\}
\]
to be linearly independent,
with $e_\beta^{(i)}\in X_{g_\beta}$ for every $\beta<\alpha$ and $i=1,2$. Since the subgroup generated by this set has cardinality strictly less than $\mathfrak c$, Lemma~\ref{Le1} yields elements
$
e_\alpha^{(1)},e_\alpha^{(2)}\in X_{g_\alpha}
$
such that $\{e_\beta^{(1)},e_\beta^{(2)}:\beta\le\alpha\}$ remains linearly independent. This completes the induction.

For $i=1,2$, set
\[
E_i=\{e_\alpha^{(i)}:\alpha<\mathfrak c\},
\qquad
Y_i=\langle E_i\rangle .
\]
Then $|E_1|=|E_2|=\mathfrak c$. Since $E_1\cup E_2$ is linearly independent, we have $Y_1\cap Y_2=\{0\}$. Finally, for every nonzero $g\in G$, say $g=g_\alpha$, we have $e_\alpha^{(i)}(g)=1$ for $i=1,2$. Hence both $Y_1$ and $Y_2$ separate points of $G$.

We now use $E_1$ and $E_2$ to define two group topologies on $G$. For each $i=1,2$, let $\tau_i$ denote the group topology on $G$ generated by $E_i$ (equivalently, by $Y_i$).
Since $Y_i$ separates points of $G$, the topology $\tau_i$ is Hausdorff.
Moreover, the map
\[
h_i:(G,\tau_i)\to 2^{E_i},
\qquad
g\mapsto\Psi_g,
\]
where $\Psi_g(\chi)=\chi(g)$ for $\chi\in E_i$, defines a topological group embedding.

Let
\[
\Pi=2^{E_1}\times2^{E_2}
\]
which we naturally identify with $2^{E_1\cup E_2}$. Consider the homomorphism
\[
h:(G,\tau)\longrightarrow\Pi,\qquad
h(g)=(h_1(g),h_2(g)),
\]
where $\tau$ is the group topology on $G$ generated by
$E_1\cup E_2$ (equivalently, by $Y_1\oplus Y_2$).
Since $E_1\cup E_2$ separates points of $G$, the map $h$ is a
topological group embedding. Set
\[
H=h(G)=\{(h_1(g),h_2(g)):g\in G\}.
\]
Thus $H$ is topologically isomorphic to $(G,\tau)$.
Since $\tau$ is weaker than the original topology of $G$, the group
$(G,\tau)$, and hence $H$, is countably compact.

Set
\[
\Pi_2=\{0\}\times 2^{E_2}\leq \Pi
\]
and let $\Sigma_2$ denote the $\Sigma$-product in $\Pi_2$, that is,
\[
\Sigma_2=\left\{(0,x)\in\Pi_2:
|\operatorname{supp}(x)|\leq\omega\right\}.
\]
Recall that a space is $\omega$-bounded if every countable subset has compact closure. It is easy to verify that $\Sigma_2$ is $\omega$-bounded.

Finally, set
\[
K=H+\Sigma_2.
\]
Since $H$ is countably compact and $\Sigma_2$ is $\omega$-bounded, $K$ is a countably compact group by \cite[Lemma~5.4]{LMT}.

Since $h_1$ is injective,
\[
H\cap\Pi_2=\{(0,0)\}.
\]
Indeed, if $(h_1(g),h_2(g))\in H\cap\Pi_2$, then $h_1(g)=0$, and hence $g=0$. Consequently,
\[
K\cap\Pi_2=\Sigma_2.
\]
The inclusion $\Sigma_2\subseteq K\cap\Pi_2$ is immediate. Conversely,
if $h+s\in K\cap\Pi_2$, where $h\in H$ and $s\in\Sigma_2$, then
\[
h=(h+s)-s\in H\cap\Pi_2.
\]
Thus $h=0$, and hence $h+s=s\in\Sigma_2$.
Since $\Pi_2$ is closed in $\Pi$ and $K\cap\Pi_2=\Sigma_2$, the subgroup $\Sigma_2$ is closed in $K$. We also have $d(\Sigma_2)=\mathfrak c$. Indeed, the finite-support subgroup of $2^{E_2}$ is dense in $\Sigma_2$ and has cardinality $\mathfrak c$, so $d(\Sigma_2)\leq\mathfrak c$. Conversely, if $D\subseteq\Sigma_2$ is dense, then
\[
E_2=\bigcup_{x\in D}\operatorname{supp}(x).
\]
Otherwise, if $\xi\in E_2$ were omitted from this union, the basic open set $\{x\in\Sigma_2:x(\xi)=1\}$ would miss $D$. Since every support is countable and $|E_2|=\mathfrak c$, it follows that $\mathfrak c\leq |D|\cdot\omega$, and hence $|D|\geq\mathfrak c$. Thus $d(\Sigma_2)=\mathfrak c$, so $\Sigma_2$ is a closed non-separable subgroup of $K$.

It remains to show that $K$ is separable. Let
\[
H_0=\{(h_1(g),h_2(g)):g\in G_0\}.
\]
Since $G_0$ is countable, so is $H_0$. 

We claim that $H_0$ is dense in $\Pi$. Since $G_0$ is dense in $G$, the restrictions to $G_0$ of the
characters in $E_1\cup E_2$ remain linearly independent. Thus, for
every finite $F\subseteq E_1\cup E_2$, the homomorphism
\[
G_0\longrightarrow 2^F,\qquad
g\longmapsto(\chi(g))_{\chi\in F},
\]
is surjective; otherwise, its image would be contained in the kernel of a nonzero linear functional on $2^F$, yielding a nontrivial linear relation among the restrictions of the elements of $F$. Hence $H_0$ meets every basic open subset of $2^{E_1\cup E_2}$, and therefore $H_0$ is dense in $\Pi$.

Since
\[
H_0\subseteq K\subseteq\Pi,
\]
it follows that $H_0$ is dense in $K$. Thus $K$ is separable.
Together with the above, this completes the $\mathbf{ZFC}$ construction and establishes Theorem \ref{ThB}.\end{proof}

\subsection{The density spectrum of closed subgroups}

We now turn to a related but more global problem, namely the determination of the \emph{density spectrum of closed subgroups} of a topological group.

\begin{definition}
For an infinite topological group $G$, the density spectrum of its closed subgroups is defined by
\[
\cd(G)=\{d(H): H\leq G \text{ is an infinite closed subgroup}\}.
\]
\end{definition}

For compact groups, this notion is closely related to the possible weights of closed subgroups. Indeed, it is well known that every infinite compact group satisfies
\[
d(G)=\log(w(G)).
\]
Consequently, for an infinite compact group $G$, the determination of $\cd(G)$ reduces to determining the possible weights of its infinite closed subgroups.

The latter problem was completely settled by Hern\'{a}ndez, Hofmann, and Morris.

\begin{theorem}\cite[Main Theorem]{HHM}\label{HHM}
Let $G$ be an infinite locally compact group. Then for every cardinal
\[
\kappa\in[\omega,w(G)],
\]
there exists a closed subgroup $H$ of $G$ such that $w(H)=\kappa$.
\end{theorem}

Combining this theorem with the identity $d(K)=\log(w(K))$ for infinite compact groups immediately yields the following description of the density spectrum.

\begin{theorem}\label{cdk}
Let $G$ be an infinite compact group. Then
\[
\cd(G)=\{\log\kappa:\kappa\in[\omega,w(G)]\}.
\]
\end{theorem}

The result above shows that, from the perspective of arbitrary closed subgroups, every possible weight between $\omega$ and $w(G)$ is realized. A natural question is whether the same phenomenon persists when one restricts attention to \emph{normal} closed subgroups.

\begin{remark}
In \cite[p.~621, Problem]{HHM}, the authors asked whether, for every compact group $G$ and every infinite cardinal $\kappa\leq w(G)$, there exists a \emph{normal} closed subgroup $N$ of $G$ such that
$
w(N)=\kappa.
$
If true, this would provide a normal-subgroup analogue of Theorem~\ref{HHM} and would strengthen Theorem~\ref{cdk} accordingly. We shall eventually show in Proposition~\ref{prop_free_prop} that the answer is negative, thereby resolving the question posed in \cite{HHM}.
\end{remark}

Before giving the example, recall that a pro-$p$ group is called \emph{free pro-$p$} if it satisfies the usual universal property in the category of pro-$p$ groups. We shall use the fact that every closed subgroup of a free pro-$p$ group is free pro-$p$; see \cite[Corollary~7.7.5]{RZ}.

\begin{proposition}\label{prop_free_prop}
There exists a profinite group without non-trivial metrizable closed normal subgroups.
\end{proposition}

\begin{proof}
Let $p$ be a prime, let $\kappa$ be an uncountable cardinal, and let $G$ be the free pro-$p$ group of rank $\kappa$. Suppose that $N$ is a non-trivial metrizable closed normal subgroup of $G$.

First suppose that $[G:N]<\infty$. Then $N$ is an open subgroup of the free pro-$p$ group $G$. By \cite[Theorem~3.6.2]{RZ}, $N$ is free pro-$p$ and, since $\kappa$ is infinite,
\[
\operatorname{rank}(N)=\operatorname{rank}(G)=\kappa.
\]

Suppose instead that $[G:N]$ is infinite. Then $N$ is free pro-$p$ by \cite[Corollary~7.7.5]{RZ}, and \cite[Proposition~8.6.3]{RZ} gives
\[
\operatorname{rank}(N)=\max\{\kappa,\omega\}=\kappa.
\]

Thus, in either case, $N$ is an infinite free pro-$p$ group of rank $\kappa$. By the equality between rank and $\delta$ established above and Lemma~\ref{delta=w}, we obtain
\[
w(N)=\delta(N)=\operatorname{rank}(N)=\kappa.
\]
On the other hand, every infinite compact metrizable group has countable weight, so $w(N)=\omega$. This contradicts $\kappa>\omega$. Thus $G$ has no non-trivial metrizable closed normal subgroup.
\end{proof}

Theorem \ref{cdk} reveals that the density spectrum of closed subgroups $\cd(G)$ for a compact group is completely determined by the logarithmic weights of its closed subgroups. However, this characterization inherently suggests that $\cd(G)$ can hardly form an unbroken interval of cardinals. For instance, under the set-theoretic assumption that $2^{\omega_1} = 2^\omega$, the cardinal $\omega_1$ cannot serve as the density of any compact group. Consequently, if a compact group $G$ has density $d(G) > \omega_1$, its density spectrum $\cd(G)$ contains a gap and is not an interval.

\subsection{Proof of Theorem \ref{ThC}}
We conclude this paper by investigating the extent to which $\cd(G)$ \emph{does} contain rich, ``almost'' unbroken segments of cardinals. Shifting our focus to general topological groups, we present a positive structural result that guarantees the existence of closed subgroups of prescribed densities. To formulate this final result, we recall the fundamental notion of topological \emph{tightness}. 

The tightness $t(X)$ of a topological space $X$ is defined as the minimal infinite cardinal $\lambda$ such that, whenever a point $x$ lies in the closure of a subset $A \subseteq X$, there exists a subset $B \subseteq A$ of cardinality $|B| \le \lambda$ such that $x \in \overline{B}$. By utilizing tightness as a topological bound, we show that any topological group contains closed subgroups realizing all regular densities between its tightness and its overall density.

\begin{proposition}\label{Prop:t}
Let $G$ be a topological group with density $\kappa > \omega$. Then for every regular cardinal $\tau$ such that $t(G) < \tau \le \kappa$, $G$ contains a closed subgroup of density $\tau$.
\end{proposition}

\begin{proof}
We recursively construct a strictly increasing sequence
$\{H_\alpha:\alpha\leq\tau\}$ of closed subgroups of $G$. Choose any
$x_0\in G$ and put $H_0=\overline{\langle x_0\rangle}$. Suppose that
$H_\beta$ and $x_{\beta+1}$ have been chosen for all $\beta<\alpha$, and set
\[
X_\alpha=\{x_0\}\cup\{x_{\beta+1}:\beta<\alpha\}.
\]
Inductively we shall have $H_\alpha=\overline{\langle X_\alpha\rangle}$.
If $\alpha<\tau$, then
\[
d(H_\alpha)\leq |\langle X_\alpha\rangle|
\leq\max\{|\alpha|,\omega\}<\tau\leq\kappa=d(G),
\]
so $H_\alpha\neq G$. Choose $x_{\alpha+1}\in G\setminus H_\alpha$ and put
\[
H_{\alpha+1}=\overline{\langle H_\alpha\cup\{x_{\alpha+1}\}\rangle}.
\]
At a non-zero limit ordinal $\alpha\leq\tau$, put
$H_\alpha=\overline{\bigcup_{\beta<\alpha}H_\beta}$. Since
$\langle X_\alpha\rangle=\bigcup_{\beta<\alpha}\langle X_\beta\rangle$
at limit stages, the identity
$H_\alpha=\overline{\langle X_\alpha\rangle}$ follows by induction.
In particular, $\langle X_\tau\rangle$ is dense in $H_\tau$ and has
cardinality at most $\tau$, so $d(H_\tau)\leq\tau$.

Suppose that $d(H_\tau)<\tau$. Choose a dense subset $A\subseteq H_\tau$ with $|A|<\tau$. Since $\bigcup_{\alpha<\tau}H_\alpha$ is dense in $H_\tau$, for each $a\in A$ there is a set
$B_a\subseteq\bigcup_{\alpha<\tau}H_\alpha$ such that
$a\in\overline{B_a}$ and $|B_a|\leq t(G)<\tau$. Put $B=\bigcup_{a\in A}B_a$. By regularity of $\tau$, $|B|<\tau$. Again using regularity, there is some $\alpha<\tau$ such that $B\subseteq H_\alpha$. Therefore
\[
A\subseteq\overline B\subseteq H_\alpha,
\]
and hence $H_\tau=\overline A\subseteq H_\alpha$, contradicting
$H_\alpha\subsetneq H_{\alpha+1}\subseteq H_\tau$. Thus $d(H_\tau)=\tau$.
\end{proof}

\begin{corollary}\label{cor:regular-spectrum}
Let $G$ be a topological group with countable tightness and let $\tau$ be a
regular infinite cardinal. If $\kappa\in\cd(G)$ for some $\kappa\geq\tau$,
then $\tau\in\cd(G)$.
\end{corollary}

\begin{proof}
Choose a closed subgroup $H$ of $G$ with $d(H)=\kappa$. If $\tau=\omega$,
the closure in $H$ of any countably infinite subgroup is an infinite
separable closed subgroup. If $\tau>\omega$, then
$t(H)\leq\omega<\tau$, and Proposition~\ref{Prop:t} applied to $H$
yields a closed subgroup of density $\tau$.
\end{proof}

For the singular case, the relevant cardinal invariant is the cofinality
of the poset of countable subsets. The following estimate is the key
additional ingredient.

\begin{lemma}\label{lem:closed-density-bound}
Let $H$ be an $\omega$-bounded topological group with $t(H)=\omega$ and
$d(H)=\lambda$. Then every closed subgroup $L$ of $H$ satisfies
\[
d(L)\leq
\max\left\{\omega,
\operatorname{cf}\bigl([\lambda]^{\leq\omega},\subseteq\bigr)
\right\}.
\]
\end{lemma}

\begin{proof}
Choose a dense set $D\subseteq H$ with $|D|=\lambda$, and let
$\mathcal F\subseteq[D]^{\leq\omega}$ be cofinal under inclusion with
$|\mathcal F|
=\operatorname{cf}([\lambda]^{\leq\omega},\subseteq)$.
For each $F\in\mathcal F$, put
$K_F=\overline{\langle F\rangle}$.

Since $\langle F\rangle$ is countable and $H$ is $\omega$-bounded,
$K_F$ is compact. Moreover, $t(K_F)\leq\omega$, and hence the compact
group $K_F$ is metrizable; see \cite{Ismail}. Thus
$L\cap K_F$ is separable.

For each $x\in L$, since $x\in\overline D$ and $t(H)=\omega$, there is a
countable set $A_x\subseteq D$ such that $x\in\overline{A_x}$. Choose
$F\in\mathcal F$ with $A_x\subseteq F$. Then $x\in K_F$, and therefore
$L=\bigcup_{F\in\mathcal F}(L\cap K_F)$.

For every $F\in\mathcal F$, choose a countable dense subset $E_F$ of
$L\cap K_F$. Then $\bigcup_{F\in\mathcal F}E_F$ is dense in $L$, and
the desired estimate follows.
\end{proof}

The estimate becomes particularly transparent when countable powers do
not increase the cardinal.

\begin{corollary}\label{cor:lambda-omega-monotone}
Let $H$ be an $\omega$-bounded topological group with $t(H)=\omega$ and
$d(H)=\lambda$. If $\lambda^\omega=\lambda$, then every closed subgroup
$L$ of $H$ satisfies $d(L)\leq\lambda$.
\end{corollary}

\begin{proof}
Since $|[\lambda]^{\leq\omega}|\leq\lambda^\omega=\lambda$, we have
$\operatorname{cf}([\lambda]^{\leq\omega},\subseteq)\leq\lambda$.
Therefore
\[
\max\left\{\omega,
\operatorname{cf}([\lambda]^{\leq\omega},\subseteq)\right\}\leq\lambda,
\]
and Lemma~\ref{lem:closed-density-bound} applies. Notice that no reverse inequality is needed; in particular, when $\lambda=\omega$ the poset $[\omega]^{\leq\omega}$ has a largest element and hence cofinality $1$.
\end{proof}

We can now extend Proposition~\ref{Prop:t} to a class of singular
cardinals.

\begin{proposition}\label{prop:singular-omega-bounded}
Let $G$ be an $\omega$-bounded topological group with $t(G)=\omega$, and
let $\tau$ be a singular cardinal. Assume that
$\operatorname{cf}([\lambda]^{\leq\omega},\subseteq)<\tau$ for every
infinite cardinal $\lambda<\tau$. If $\kappa\in\cd(G)$ for some
$\kappa\geq\tau$, then $\tau\in\cd(G)$.
\end{proposition}

\begin{proof}
Choose a closed subgroup $C$ of $G$ with $d(C)=\kappa$. Let
$\{\tau_i:i<\operatorname{cf}(\tau)\}$ be an increasing cofinal family
of uncountable regular cardinals below $\tau$. By
Proposition~\ref{Prop:t}, for each
$i<\operatorname{cf}(\tau)$ there is a closed subgroup $L_i$ of $C$
such that $d(L_i)=\tau_i$.

Choose a dense set $D_i\subseteq L_i$ with $|D_i|=\tau_i$ and put
\[
H=
\overline{\left\langle
\bigcup_{i<\operatorname{cf}(\tau)}D_i
\right\rangle}.
\]
Then $d(H)\leq\tau$. Moreover, $D_i\subseteq H$ and $H$ is closed in
$C$, so $L_i=\overline{D_i}\subseteq H$. Hence each $L_i$ is a
closed subgroup of $H$.

Suppose that $d(H)=\lambda<\tau$. By
Lemma~\ref{lem:closed-density-bound}, for every
$i<\operatorname{cf}(\tau)$,
\[
\tau_i=d(L_i)\leq
\max\left\{\omega,
\operatorname{cf}\bigl([\lambda]^{\leq\omega},\subseteq\bigr)
\right\}<\tau.
\]
The cardinal on the right is independent of $i$, contradicting the
cofinality of $\{\tau_i:i<\operatorname{cf}(\tau)\}$ in $\tau$.
Thus $d(H)=\tau$, and so $\tau\in\cd(G)$.
\end{proof}

A convenient sufficient condition is obtained from ordinary cardinal
exponentiation.

\begin{corollary}\label{cor:singular-power-condition}
Let $G$ be an $\omega$-bounded topological group with countable tightness,
and let $\tau$ be singular. Assume that $\lambda^\omega<\tau$ for every
infinite cardinal $\lambda<\tau$. If $\kappa\in\cd(G)$ for some
$\kappa\geq\tau$, then $\tau\in\cd(G)$.
\end{corollary}

\begin{proof}
For every infinite $\lambda<\tau$, we have
$\operatorname{cf}([\lambda]^{\leq\omega},\subseteq)
\leq |[\lambda]^{\leq\omega}|
\leq\lambda^\omega<\tau$.
Apply Proposition~\ref{prop:singular-omega-bounded}.
\end{proof}

\begin{corollary}\label{cor:singular-strong-limit}
Let $G$ be an $\omega$-bounded topological group with countable tightness,
and let $\tau$ be a singular strong limit cardinal. If
$\kappa\in\cd(G)$ for some $\kappa\geq\tau$, then $\tau\in\cd(G)$.
\end{corollary}

\begin{proof}
For every infinite $\lambda<\tau$,
$\lambda^\omega\leq(2^\lambda)^\omega=2^\lambda<\tau$.
Hence Corollary~\ref{cor:singular-power-condition} applies.
\end{proof}

The first singular cardinal can be handled in $\mathbf{ZFC}$ without any
strong limit assumption.

\begin{corollary}\label{cor:aleph-omega}
Let $G$ be an $\omega$-bounded topological group with countable tightness.
If $\kappa\in\cd(G)$ for some $\kappa\geq\aleph_\omega$, then
$\aleph_\omega\in\cd(G)$.
\end{corollary}

\begin{proof}
It is standard that
$\operatorname{cf}([\aleph_n]^{\leq\omega},\subseteq)=\aleph_n$ for
every $1\leq n<\omega$. For completeness, this follows by induction.
The case $n=1$ uses the fact that every countable subset of $\omega_1$
is bounded. If the assertion holds for $n$, then every countable subset
of $\aleph_{n+1}$ is contained in some $\alpha<\aleph_{n+1}$, while
$|\alpha|\leq\aleph_n$; taking the union of cofinal families for these
initial segments gives the upper bound $\aleph_{n+1}$, and the reverse
inequality is immediate.

Therefore
$\operatorname{cf}([\lambda]^{\leq\omega},\subseteq)<\aleph_\omega$
for every infinite cardinal $\lambda<\aleph_\omega$. Apply
Proposition~\ref{prop:singular-omega-bounded}.
\end{proof}

The preceding criterion also interacts naturally with pcf theory.
We briefly clarify the two forms of Shelah's Strong Hypothesis that appear
below. The usual, global form of Shelah's Strong Hypothesis
$\mathbf{SSH}$ asserts that $pp(\mu)=\mu^+$ for every singular cardinal
$\mu$, where $pp(\mu)$ is the pseudo-power; see \cite{ShelahCA}. In terms of countable covering numbers, one of
its equivalent consequences is that
$\operatorname{cf}([\mu]^{\leq\omega},\subseteq)=\mu^+$ whenever
$\mu$ is uncountable and $\cof(\mu)=\omega$; see also
\cite{EFHST}.

For a fixed singular cardinal $\tau$, however, our argument requires only
the corresponding local condition below $\tau$, namely
$\operatorname{cf}([\mu]^{\leq\omega},\subseteq)=\mu^+$ for every
uncountable $\mu<\tau$ with $\operatorname{cf}(\mu)=\omega$.
We shall refer to this as the \emph{local $\mathbf{SSH}$ condition below
$\tau$}. This is formally weaker than the global hypothesis used later:
it is restricted to cardinals below the particular $\tau$ under
consideration and only to those of countable cofinality. Thus global
$\mathbf{SSH}$ implies the local condition below every singular
$\tau$, while the latter is all that is needed to obtain the corresponding
conclusion for that particular $\tau$.

\begin{corollary}\label{cor:local-ssh}
Let $G$ be an $\omega$-bounded topological group with countable tightness,
and let $\tau$ be singular. Assume local $\mathbf{SSH}$ condition below
$\tau$. If $\kappa\in\cd(G)$ for some
$\kappa\geq\tau$, then $\tau\in\cd(G)$.
\end{corollary}

\begin{proof}
Write
\[
q(\lambda)=\operatorname{cf}([\lambda]^{\leq\omega},\subseteq).
\]
We prove by induction on the infinite cardinals $\lambda<\tau$ the following two assertions simultaneously:
$q(\lambda)\leq\lambda^+$, and, whenever $\operatorname{cf}(\lambda)>\omega$, in fact $q(\lambda)\leq\lambda$.

For $\lambda=\omega$, we have $q(\omega)=1$. If $\lambda>\omega$ and $\operatorname{cf}(\lambda)=\omega$, the local $\mathbf{SSH}$ hypothesis gives $q(\lambda)=\lambda^+$.

Suppose now that $\operatorname{cf}(\lambda)>\omega$. Every countable subset of $\lambda$ is bounded in $\lambda$. For each $\alpha<\lambda$, let $\mu_\alpha=\max\{\omega,|\alpha|\}<\lambda$. By the induction hypothesis there is a cofinal family in $[\alpha]^{\leq\omega}$ of cardinality at most $\mu_\alpha^+\leq\lambda$ (transporting a cofinal family along a bijection between $\alpha$ and $|\alpha|$ when necessary). The union of these families over all $\alpha<\lambda$ has cardinality at most $\lambda$ and is cofinal in $[\lambda]^{\leq\omega}$. Hence $q(\lambda)\leq\lambda$.

Thus, for every infinite $\lambda<\tau$, we have $q(\lambda)\leq\lambda^+<\tau$; the last inequality holds because $\tau$ is singular and therefore is not a successor cardinal. Proposition~\ref{prop:singular-omega-bounded} now applies.
\end{proof}

Under the full Strong Hypothesis, the $\omega$-bounded case of the
interval problem has a complete positive solution.

\begin{proof}[Proof of Theorem \ref{ThC}]
Let $\kappa\in\cd(G)$ and let $\omega\leq\tau\leq\kappa$. If $\tau$ is
regular, Corollary~\ref{cor:regular-spectrum} gives
$\tau\in\cd(G)$. If $\tau$ is singular, by Corollary \ref{cor:local-ssh} again $\tau\in\cd(G)$.
Therefore $\cd(G)$ is an interval. 
\end{proof}

Thus, in $\mathbf{ZFC}$ the possible obstruction for an
$\omega$-bounded group of countable tightness is concentrated at
singular cardinals $\tau$ for which
$\operatorname{cf}([\lambda]^{\leq\omega},\subseteq)$ can reach
$\tau$ for some $\lambda<\tau$. The preceding corollary shows that this
obstruction disappears under $\mathbf{SSH}$. For general infinite countably
compact groups, however, the compact-metrizable hull argument used in
Lemma~\ref{lem:closed-density-bound} is no longer available. This leaves
the following question open.

\begin{question}
Let $G$ be an infinite countably compact topological group with countable tightness.
Is $\cd(G)$ necessarily an interval of cardinals? In $\mathbf{ZFC}$,
is the same true for every infinite $\omega$-bounded group of countable tightness?
\end{question}
\section*{Acknowledgements}

The authors would like to express sincere gratitude to Wei He, Arkady G. Leiderman, and Mikhail G. Tkachenko for their valuable comments and suggestions on this manuscript. Special thanks are due to Professor Leiderman for suggesting the reference to the construction in \cite{HMRS}, which served as the key inspiration for resolving the open problem posed in \cite{LMT}. The authors are also grateful to Dr. Jiaheng Zhuang for his insightful discussions concerning this construction.

Work on the second part of this paper began during the second author's visit to the Institute of Mathematics at Nanjing Normal University in 2025. She gratefully acknowledges the hospitality of the Institute and the financial support of the National Natural Science Foundation of China (grants Nos.~12301089 and 12271258).

\subsection*{Declaration on the use of generative AI}
Theorem \ref{ThA} was first announced at ATTI 2025 in Udine. During a subsequent verification of the manuscript with the assistance of GPT-5.6 Sol, a serious gap was identified in the original proof of the profinite case. With assistance from the same system, the authors developed and checked a substantially revised argument that repairs this gap. Generative AI tools, primarily GPT-5.6 Sol, were also used throughout the preparation of the manuscript for language polishing, consistency checks, and verification of arguments. Some expository passages introducing standard concepts were initially drafted with AI assistance. All AI-assisted mathematical arguments, statements, references, and expository material were independently checked and verified by the authors, who take full responsibility for the content of the paper.

\end{document}